\documentclass[11pt,a4paper]{article}
\usepackage{hyperref}
\usepackage{tikz-cd}
\usepackage[utf8]{inputenc}
\usepackage{orcidlink}
\usepackage{graphicx}  \usepackage{mathptmx}
\usepackage{enumitem}
\usepackage{amsmath, amssymb}
\numberwithin{equation}{section}
\usepackage{gensymb}

\newcommand{\id}[1]{\ensuremath{\mathrm{id}}}

\newcommand{\beq}{\begin{equation}}
\newcommand{\eeq}{\end{equation}} 
\newcommand{\bea}{\begin{eqnarray}}
\newcommand{\eea}{\end{eqnarray}}

 \newcommand{\phv}{\varphi}
\newcommand{\ch}{\ch}

\newcommand{\N}{{\mathbb N}} 
 
  \makeatletter
 
\makeatletter
\def\moverlay{\mathpalette\mov@rlay}
\def\mov@rlay#1#2{\leavevmode\vtop{%
   \baselineskip\z@skip \lineskiplimit-\maxdimen
   \ialign{\hfil$\m@th#1##$\hfil\cr#2\crcr}}}
\newcommand{\charfusion}[3][\mathord]{
    #1{\ifx#1\mathop\vphantom{#2}\fi
        \mathpalette\mov@rlay{#2\cr#3}
      }
    \ifx#1\mathop\expandafter\displaylimits\fi}
\makeatother

\newcommand{\W}{Wittgenstein}

\topmargin = - 1 cm \textheight = 23 cm \textwidth = 15 cm
\usepackage[symbol]{footmisc}
\renewcommand{\thefootnote}{\fnsymbol{footnote}}
\begin{document}
\pagenumbering{arabic} \setlength{\unitlength}{1cm}\cleardoublepage
\date\nodate
\begin{center}
\begin{Huge}
{\bf Is mathematics like a game?}
\end{Huge}
\bigskip\bigskip

\begin{Large}
 Klaas Landsman\orcidlink{0000-0003-2651-2613} and Kirti Singh\vspace{5mm}
 \end{Large}
 
 \begin{large}
Institute for Mathematics, Astrophysics, and  Particle Physics\\ \vspace{1mm}
 Radboud Center for Natural Philosophy  \\ \vspace{1mm}
Radboud University, Nijmegen, The Netherlands\\ \vspace{1mm}
\texttt{landsman@math.ru.nl}\mbox{     }  \texttt{kirtiks36@gmail.com}
\end{large}\bigskip

 \begin{abstract} 
\noindent 
We re-examine the old question to what extent mathematics may be compared with a game. Mainly inspired by Hilbert and \W, our answer is that mathematics is something like  a ``rhododendron of language games'', where the rules are inferential. The pure side of mathematics is essentially formalist, where we propose that truth is not carried by theorems  corresponding to whatever independent reality and 
  \emph{arrived at} through proof, but is  \emph{defined by}  correctness of rule-following (and as such is  objective given these rules).  G\"{o}del's  theorems, which are often seen as a threat to formalist philosophies of mathematics, actually strengthen our concept of truth. 
The applied side of mathematics arises from two practices: first, 
 the dual nature of axiomatization as \emph{taking} from heuristic practices like physics and informal mathematics whilst \emph{giving} proofs and logical analysis; and second, the ability of  using the inferential role of theorems to make
 ``surrogative'' inferences about natural phenomena. Our framework 
 is pluralist,  combining various  (non-referential) philosophies of mathematics. \end{abstract}\end{center}
\tableofcontents

\medskip

\noindent\textbf{Statements and Declarations:} \emph{The authors have no competing interests}.
\medskip

\noindent \textbf{Acknowledgement:} This paper originated in a B.Sc.\ thesis of the second author, which has been
 greatly expanded by the first. The authors thank Simon Friederich, Barteld Kooi,  Tushar Menon, Felix M\"{u}hlh\"{o}lzer, Fred Muller, Tomasz Placek, Chris Scambler, Robert Thomas,  Freek Wiedijk, and five anonymous referees (for various journals!) for comments and  help. The first author also thanks All Souls College, Oxford, where the paper was finished in June 2025, for its hospitality.
\thispagestyle{empty}
\newpage
\renewcommand{\thefootnote}{\arabic{footnote}}
 \setcounter{footnote}{0}
\section{Introduction}\label{Intro}
The  aim of this paper is to re-examine the old question to what extent mathematics may be compared with a game (like chess).\footnote{The review by Epple (1994) provides an excellent historical and philosophical introduction to this question; see also Detlefsen (2005) and Weir (2022). We will not discuss the the  mathematics of games, cf.\ du Sautoy (2023).} It is no accident that  the idea of such a comparison originated in the late nineteenth century, since that was the time of the ``modernist transformation'' in which mathematics  lost its connection with physical reality and visualizability, which were replaced by abstraction and rigorous proof.\footnote{See, for example, Mehrtens (1990) and  Gray (2008).  \label{Grayfn}}
 Indeed, serious analysis of the analogy between mathematics and games like chess started with Frege's criticisms of Thomae (1898), Heine (1872) and Illigens (1893).\footnote{See Frege (1903), \S\S 86--137, translated in Geach and Black (1960). } 
Frege  rejected the analogy;  he was primarily unable to comprehend how a mere game could describe any ``thought'',\footnote{ This was the central ingredient of his philosophy of both mathematics and language. The word `thought' (\emph{Gedanke})  is ambiguous in both English and German. Frege (1892)  clarifies that for him, thoughts do not refer to the subjective act of thinking, but to the objective content thereof, and thus are possible truth carriers. Frege (1918) moves to an uncompromising platonism, e.g.,  `Without wishing to give a definition, I call a thought something for which the question of truth arises. (\ldots) A thought is something immaterial and everything material and perceptible is excluded from this sphere of that for which the question of truth arises.' (Frege, 1918, p.\ 292). He even uses the (no longer usable) term `dritte Reich' (best translated as: third Realm) as the sphere beyond the material and perceptible. 
 } and his secondary, more specific reasons for  disapproving of the analogy  reduce to this idiosyncratic inability (see \S\ref{FHW} below).  
 Frege's discussion remains of considerable interest for the general philosophy of mathematics, whose  questions concern:\footnote{For a somewhat different list see Linnebo (2017), \S 1.1.  See also Tait (2001).}
\begin{enumerate}
\item \emph{Ontology:} What is mathematics \emph{about}?\footnote{Anticipating a \W ian setting, a more precise form of this question would be:  \emph{are  mathematical objects given in advance?}, where 
we adopt the following definition:  `An object is \emph{given in advance} iff the criteria of identity for the object which the language game is about are \emph{not} completely stated or presented by the language game itself; and it is  \emph{not given in advance} iff the criteria of identity for the object are completely stated or presented in the language game -- [so] that this identity is given by the language game alone and by nothing else.'  (M\"{u}hlh\"{o}lzer, 2012, p.\ 114).} What (and where) are mathematical objects?
\item \emph{Truth:}  What is the nature of mathematical truth?  
\item \emph{Epistemology:} How  can we know about mathematics? 
 Do we \emph{discover} or \emph{invent} it?
 \item \emph{Applicability:}  What makes applied mathematics  possible?
\end{enumerate}
We cannot even begin to summarize the huge literature about these questions, which go back to Plato and Aristotle. But there is one position--incorporating all kinds of platonism and other forms of realism--we would like to highlight as a foil, since it seems to be very widely shared among both ``working mathematicians'' and philosophers. Here is a typical expressions of it by Hardy:\footnote{See also Hardy (1940), pp.\ 63--64, or more recently: `The assumption that an arithmetical statement
is true is not an assumption about what can be proved in any formal system, or about what can be ``seen to be true,'' and nor is it an assumption
presupposing any dubious metaphysics. Rather, the assumption that Goldbach's conjecture is true is exactly equivalent to the assumption that every even number greater than 2 is the sum of two primes. Similarly, the assumption that the twin prime conjecture is true means no more and no less than the assumption that there are infinitely many primes $p$ such that
$p+2$ is also a prime, and so on. In other words ``the twin prime conjecture
is true'' is simply another way of saying exactly what the twin prime conjecture says. It is a mathematical statement, not a statement about what can be known or proved, or about any relation between language and a
mathematical reality.' (Franz\'{e}n, 2005, p.\ 30)}
\begin{quote}\begin{small}
`It seems to me that no philosophy can possibly be sympathetic to a mathematician which does not admit, in one manner or another, the immutable and unconditional validity of mathematical truth. Mathematical theorems are true or false; their truth or falsity is absolute and independent of our knowledge of them. In some sense, mathematical truth is part of objective reality. ``Any number is the sum of 4 squares''; ``any number is the sum of 3 squares''; ``any even number is the sum of 2 primes''. These are not convenient working hypotheses, or half-truths about the Absolute, or collections of marks on paper, or classes of noises summarising reactions of laryngeal glands. They are, in one sense or another, however elusive and sophisticated that sense may be, theorems concerning reality, of which the first is true, the second is false, and the third is either true or false, though which we do not know. They are not creations of our minds; Lagrange discovered the first in 1774; when he discovered it he discovered some- thing; and to that something Lagrange, and the year 1774, are equally indifferent.' (Hardy, 1929, p.\ 4).
\end{small}
\end{quote}
This belief that any proposition about at least natural numbers is either true or false \emph{full stop}, i.e., irrespective of axioms, deductions, etc.,
 has been called \emph{arithmetical determinacy} (Warren, 2020) or \emph{truth-completeness of arithmetic} (Paseau and Pregel, 2023). The latter also call it  a `\emph{fundamental commitment of mathematics}'. It is debatable whether this commitment applies just to arithmetic or to further or even all parts of mathematics (as suggested by Hardy's introductory prose but not by his examples); but this difference is irrelevant to us, since we reject it even in its--most basic and most convincing-- arithmetical form. 
As we will explain, this rejection and the alternative view we replace it by are a direct consequence of our analysis of the title question, which we initially follow up through Frege and his great adversaries Hilbert and \W.

Via this route we arrive at an answer to the effect that although comparing  mathematics to a \emph{game} is far too simple, 
it may  be favourably compared with a certain combination of \emph{language games} whose rules are inferential and whose interlocking is well described by a \emph{rhododendron}, having both multiple roots and branches.\footnote{See \S\ref{MMLG} for the connection of this idea with \W's famous `motley' of language games.} 
At its roots one finds foundational theories or overall formal and logical frameworks for mathematics like ZFC set theory, \emph{as well as its serious competitors}.\footnote{Their axioms are what Feferman (1999) calls \emph{foundational}. The others are \emph{structural}. See also Schlimm (2013).}
Each of these roots branches out into one or even various forms of metamathematics, as well as into (typically interwoven) individual areas of mathematics  (like number theory or group theory), which in turn brach out into increasingly specialized areas thereof (like algebraic number theory or Lie groups). 
 Within all of these areas (as well as within the foundational theories at the bottom) one has a language 
game of pure mathematics with its formalized proofs; in many of them, one has a second language game of applied mathematics (including mathematical physics). 
 Our picture is pluralist and provides a coat rack onto which various (originally partly hostile) philosophies of mathematics may be attached and may even peacefully support each other, e.g.:
 \begin{itemize}
\item \emph{Formalism}, and its close relatives  \emph{deductivism} and \emph{coventionalism} as the language game governing pure mathematics grounded in proofs and defining the meaning of truth.
 \item \emph{Intuitionism} comes in twice: (i) as one of the possible language games giving a foundation of pure mathematics; and (ii): even if the latter is classical, as a possible logic of the metamathematics used to analyze the framework (and similarly for \emph{finitism}, as Hilbert  tried).
\item 
\emph{Inferentialism} also enters twice: first, in  pure mathematics as the source of its conventionalism, and second, in applied mathematics
through its use in the concept of \emph{surrogative inference} in relating mathematical models to empirical phenomena
(see below).
\end{itemize}
One advantage of our approach is that it describes both the \emph{practice} and the  \emph{results} of mathematics.
For example, each step in a proof is seen as a \emph{move} in a specific language game, which is a practice: the final \emph{score} of the game is  the result. For this reason we also incorporate the
 \emph{philosophy of mathematical practice}, which aligns with our aim (which was also \W's) to describe mathematics \emph{as it is}, seeing it as a human practice grounded in history.\footnote{For introductions to the philosophy of mathematical practice we refer to
 e.g.\ Mancosu (2008),  Hamami \& Morris (2020). P\'{e}rez-Escobar (2022) argues that this philosophy resonates well with  the  later Wittgenstein.} More generally,
rather than contradicting each other all these philosophies in fact complement and reinforce each other.\footnote{In his dispute with Brouwer, Hilbert  came to see this at least for the first and the second points (Mancosu, 1998).} 
 
 As an introduction to our proposal, let us summarize our answers to the four questions above.
 \begin{enumerate}
\item \emph{Ontology}.
We follow \W\ in warning against the confusion that mathematical expressions are (primarily) referential (see \S\ref{FHW}). Numerous problems are created by assuming the ethereal ``existence'' of mathematical objects (not just in the Platonic sense). The things mathematics \emph{talks} about are similar to chess pieces, whose physical or spiritual embodiment is  irrelevant. Yet
mathematics \emph{is} not about these pieces themselves in whatever incarnation. It is about the rules they are subject to, and their consequences.  As we see it, this makes mathematics  intersubjective on the verge of objectivity in the following sense:\footnote{Rota, a mathematician, here summarizes a key element of
Husserl's  philosophy of mathematics as he saw it shortly before his death in 1999, after 40 years of study. An important  primary source is Husserl (1954).  
 From the vast secondary literature we just refer to Hacking (2009) for the work just cited, and to
  Hartimo (2021) more generally.   }
    \begin{quote}
\begin{small} 
The constitutive property of mathematical items is not existence, but identity. (\ldots) It is painful to abandon the age-old prejudice that
identity must presuppose existence. The permanence of the identity of a mathematical
item through space and history, and across civilizations, is an extraordinary
phenomenon for which there is no easy explanation, and which is shared by few
objects of the world.  (Rota, 2000, pp.\ 93) 
\end{small}
  \end{quote}
  Husserl (1954)  explained an interesting aspect of this permanence, which we endorse.  Historically, mathematics  originated in experience and applications, but subsequently underwent a process of `idealization', in that for example mathematicians idealize different drawings of a circle into identical circles.
After this first step of `idealization' (which is relevant in the specific context of Euclidean geometry that he discusses but can be replaced by formalizing any piece of would-be mathematics),  the `objectification' or  `permanence'  (\emph{Immerfort-Sein}) inherent in the concept of `identity' takes place within humanity seen as an `emphatic and linguistic community' (\emph{Einf\"{u}hlungs- und Sprachgemeinschaft}) or `communication community' (\emph{Mitteilungsgemeinschaft}). Permanence then arises via written or other communication, followed by `reactivation' (\emph{Reaktivierung}). This seems a valid description of the sense in which both pure mathematics and (tournament) chess are shared.\footnote{Dawkins's well-known concept of a \emph{meme} also seems to describe the nature of mathematics as we just described.}
  \item \emph{Truth}. Theorems lack essential properties of things that have a truth value. 
 What may be true (or a matter of fact) is the claim \emph{that some sentence is a theorem within a specific formal system}:  it is not the \emph{content of a theorem} that is true but its ``theoremhood'' (i.e., that the result follows from the premisses according to the rules; that the game has correctly been played). 
 Although this kind of truth is remote from the Platonist one, it shares the advantage  that we should all agree about it: even Brouwer should admit that Hilbert's theorems are correct \emph{by his own standards}, and \emph{vice versa}. 
 We extensively argue this point in \S\ref{truth}.
  \item \emph{Epistemology}. We \emph{know} about mathematical items  because we \emph{invented} them. 
   If mathematics is primarily seen as a human practice, the epistemological question hardly arises. For example,  if one believes with  Russell that (mathematical) logic reveals the logical structure of reality, the question is how we know this structure. For us,  any kind of logic is an inferential language game invented by logicians on the basis of studying actual and practical linguistic and mathematical reasoning, which may subsequently be studied by itself, or may be held against language or mathematics as an object of comparison, perhaps even in a normative way (as in prescribing rules of proof).\footnote{We follow late \W\ here. See for example \emph{Philosophical Investigations}, \S\S130--131,
   as well as  Bangu (2018), Kuusela (2019), chapter 6, and Peregrin (2019) for expositions. 
   More generally, as we shall see in \S\ref{FHW},  mathematical theorems are held against  practices like someone computing $25\times 25$ or counting apples, but this time as \emph{norms}. }  The epistemology of this is unproblematic.  
 \item \emph{Applicability.} Though sometimes ignored, this is one of the ``hard problems'' of philosophy:
 \begin{quote}
\begin{small} To an unappreciated degree, the history of Western Philosophy is the history of attempts to understand why mathematics is applicable to Nature, despite apparently good reasons to believe that it should not be. (Steiner, 2005, p.\ 625)
\end{small}
\end{quote}For us, the main problem is to match our non-referential account of pure mathematics with its apparent representational role in describing the natural  world. Our answer will be presented in historical and philosophical detail in  successor papers (Landsman, 2025ab).\footnote{Here is already a sketch:
 We just saw that logical language games are both extracted from natural language and held against it. Similarly, in Hilbert-style mathematical physics axioms are inspired by (heuristic) theoretical physics so as to define 
 theories of of pure mathematics (again seen as inferential language games), which are subsequently 
 held against the relevant natural phenomena as  yardsticks or objects of comparison.  But how is this comparison or measurement  actually made? This  is done via what is called \emph{surrogative inference}, in which inferences made from a mathematical model mirror inferences made about natural (or artificial) phenomena  (Su\'{a}rez, 2024). This obviously squares  with our inferential view of pure mathematics, and (perhaps less obviously) also matches our anti-realism about (pure) mathematics with a corresponding  empiricist philosophy of science as developed notably by van Fraassen (2008).}
 \end{enumerate}
 
We elaborate our proposal in three sections. In \S\ref{FHW} we summarize Frege's objections to the analogy between mathematics and games, and introduce our main protagonists Hilbert and \W\ through their replies to Frege, followed by  a summary of (late) \W, including the Brandom\-ian turn to
\emph{inferential} rule-following, and a review of Hilbert's views on the foundations of mathematics in so far as these are relevant.
\S\ref{MMLG}  presents our picture of mathematics as a  ``rhododendron of language games'', followed by an analysis of  truth in \S\ref{truth}. We conclude the main body of the paper  in \S\ref{Conc}, followed by an appendix that relates our proposal to four related philosophies of mathematics, viz.\  formalism, pluralism, conventionalism, and deductivism. 
\section{From Frege to Hilbert and \W}\label{FHW}
Our proposal that mathematics is like a rhododendron of inferential language games originated in an analysis of 
Frege's arguments against an analogy between mathematics and chess, which we  compared with the pertinent positions of Hilbert and \W.  Predicated on his idea that mathematics is about ``thoughts'' which games are allegedly too poor to carry, Frege argued that:
\begin{itemize}
\item Mathematics is \emph{meaningful}, since it refers to  thoughts.
 But games are meaningless.
\item Thus the rules of mathematics originate in reality, whereas for games they are arbitrary.
\item Grounded in reality, there is \emph{truth} in mathematical theorems, which games lack.
\item  On Frege's conceptualization, logical inference must take us from truth to truth. Without truth (which games allegedly lack) there is no logical deduction and hence no mathematics.
\item The applicability of mathematics would be incomprehensible if it were merely a game, and also leads to irresolvable ambiguities between the formal and the applied sides.
\item Even if mathematics initially were just a game, it
  also incorporates \emph{the theory of the game} instead of only \emph{being the game}. Similarly, there are theorems \emph{about} chess.\footnote{For example, a Bishop cannot move to a square of a different color: this is not a rule but a  consequence of the rules (albeit a trivial one). A deeper example is the fact that chess games must end after finitely many moves.} This is something an allegedly meaningless  game by itself (i.e.\ mathematics) could not accomplish. 
\end{itemize}
The following exchange arguably contains the essence of the debate between Frege and Thomae:\footnote{This debate (which  started in a friendly way but eventually degenerated into an acrimonious personal polemic) consists of Frege (1903), \S\S 86--103, Thomae (1906ab, 1908),  and Frege (1906, 1908ab), so that our main text only gives a very small (literally quoted) excerpt.  See also Frege (1899) in response to Schubert, written in a similar style. }
\begin{quote}
\begin{small}
Anyone who wants to ground arithmetic in a formal theory of numbers, a theory that does not ask what numbers are and what they are supposed to do, but rather asks what we need from numbers in arithmetic, will want to look at another example of purely formal creation of the human mind. I thought I had found such an example in the game of chess. The chess pieces are symbols that have no other content in the game than what is assigned to them by the rules of the game. Saying that the signs are empty may lead to misunderstandings in the absence of any goodwill to understand. So I also believed that I could view the numbers in arithmetic, seen as a game of computation, as symbols that have no other content in the game than what is assigned to them by the rules of the game or the calculation. The  system of symbols of the arithmetic game is made up of the symbols 0 1 2 3 4 5 6 7 8 9 in the usual known manner.
(Thomae, 1906a, pp.\ 434--435)
 
Mr.\ Thomae writes: `The symbol system of the arithmetic game is made up of the characters 0 1 2 3 4 5 6 7 8 9 in a known manner.' If he had simply said that the arithmetic game had  those numbers as game objects, we would be satisfied. But now he seems to want to say that the game objects are made from these numbers, and that in a known way. How should we know the matter since we first want to get to know the arithmetic game? Here Mr.\ Thomae makes the recurring mistake of assuming that what he wants to lay the foundation for, is already known.
 (Frege, 1908a, p.\ 52)\end{small}
\end{quote}
About a decade earlier, Frege had had a more balanced exchange with Hilbert on similar themes. Their correspondence between 1895 and 1903 is a jewel in the history and philosophy of mathematics, and although
 many readers will be familiar with it we now quote two passages whose theme is the same as in the exchange just cited (although chess is not mentioned explicitly):\footnote{The letters may be found in the original German in Gabriel, Kambartel, and Thiel (1980), with English translations in  Gabriel \emph{et al.} (1980).
See also   Blanchette (2018) and Rohr  (2023). }
 \begin{quote}
\begin{small}
In my opinion, a concept can be fixed logically only by its relations to other
concepts. These relations, formulated in certain statements, I call axioms,
thus arriving at the view that axioms (perhaps together with propositions
assigning names to concepts) are the definitions of the concepts.
 I did not think of this view because I had nothing better to do, but I found myself forced into it by the requirements of strictness in logical inference and in the logical construction of a theory. I have become convinced that the more subtle parts of mathematics and the natural sciences can be treated with certainty only in this way; otherwise one is only going around in a circle.
(Hilbert to Frege, 22 September 1900)
 \end{small}
 \end{quote}
Here, following his famous \emph{Grundlagen der Geometrie} from 1899,  Hilbert replaced the traditional understanding of mathematical objects (according to which they are always defined explicitly prior to appearing in axioms) by what was later called \emph{implicit definition}, according to which such objects are defined by the axiom systems in which they occur.\footnote{See Peckhaus (1996), Pollard (2010), 
Schlimmm (2011, 2013), Giovannini \& Schiemer (2021), Biagioli (2024), and 
Sereni (2024) for various aspects of the history of implicit definitions.  Briefly, there was a French line of development starting with Gergonne and summiting in Poincar\'{e}, a German one in which Pasch and to some extent Schr\"{o}der were key predecessors of Hilbert,  and an Italian line involving  Burali-Forti, Peano, and Enriques (who introduced the term `implicit definition' in the axiomatic context, which is different from Gergonne's who also used this term).  
 Giovannini \& Schiemer (2021)  call implicit definitions \emph{structural}, since the words `implicit' and `explicit' are  adjectives for certain technical definitions in logic (which  Beth's definability theorem  identifies).
  In fact, even an implicit definition \`{a} la Hilbert and Enriques may be seen as an explicit definition of an $n$-place relation, where $n$ is the number of symbols implicitly defined by the axioms (e.g.\ Blanchette, 2018, \S2.1).  
 Poincar\'{e}'s conventionalism is also based on implicit definitions (which he called `definitions in disguise' given by axioms), and despite some differences his debate with Russell  mirrored the one between Hilbert and Frege. See e.g.\  Ben-Menahem (2006).  \label{IDfn}}  
   Modern mathematics would be unthinkable without this idea, which in particular put an end to the  struggles by Cantor, Frege, and others to define sets more explicitly--a struggle that also failed for the concept of number, even for the number one, as Frege fatefully found out.  In our view, this should also have put an end to  arithmetical determinacy (cf.\ the Introduction), but it hasn't; we take this up in \S\ref{truth}. 

 Yet Hilbert was no philosopher;\footnote{This is not to say that Hilbert was a novice to philosophy; he had clearly read Kant and knew for example Husserl in person. Hilbert's student Weyl went well beyond this, but 
 as explained by Toader (2011), he never overcame the tension between: (i) the pull of Hilbert's formalism, which Weyl saw as necessary for both objectivity (in the sense of mind-independence) and free concept formation by `symbolic construction' of the kind needed for theoretical physics, yet at the cost of intelligibility;  (ii) Husserl's phenomenology, which required contentual reasoning and concept formation by abstraction from immediate experience; and (iii) Brouwer's intuitionism, which emphasized individual mathematical understanding at the cost of objectivity and formalism. 
See also Da Silva (2017).}  
 in that direction we turn to \W.\footnote{For comparisons of Hilbert and \W\ see e.g.\ Muller (2004),  M\"{u}hlh\"{o}lzer (2006, 2008, 2010, 2012), and Friederich (2011, 2014). The closest links are between Hilbert's views of mathematics around 1900 as reviewed above and those of ``middle'' \W\ (i.e., 1929--1936), where they both opposed Frege such that Hilbert's mathematics matched \W's philosophy.  In his later period \W\ became quite critical of Hilbert's metamathematics.
 }  In this section we merely discuss  \W's most acute comments on Frege and 
Hardy:
  \begin{quote}
\begin{small}
Consider [Hardy (1929)] and his remark that `to mathematical propositions there corresponds---in some sense, however sophisticated--a reality'. (\ldots)
We have here a thing which constantly happens. The words in our language have all sorts of uses; some very ordinary uses which come into one's mind immediately, and then again they have uses that are more and more remote. For instance, if I say the word `picture', you would think first and foremost of something drawn and painted and, say, hung up on the wall. You would not think of Mercator's projection of the globe; still less of the sense in which a man's handwriting is a picture of his character. A word has one or more nuclei of uses which come into every body's mind first; so that if one says so-and-so is also a picture---a map or \emph{Darstellung}
in mathematics---in this lies a comparison, as it were, ``Look at this as a continuation of that.'' So if you forget where the expression ``a reality corresponds to'' is really at home--- What is ``reality''? We think of ``reality'' of something we can \emph{point} to. It is \emph{this, that}.
(\W, \emph{Lectures on the Foundations of Mathematics}, pp.\ 239--240).
 \end{small}
 \end{quote}
 This is one of the many expressions of \W's non-referential approach to the philosophy of mathematics,\footnote{
 The philosophy of mathematics occupied \W\ throughout his career but it was never completed or even prepared for publication by \W\ himself. The primary sources are the posthumous  \emph{Bemerkungen \"{u}ber die Grundlagen der Mathematik} (BGM), written 1937--1944 (Wittgenstein, 1969) and the \emph{Lectures on the Foundations of Mathematics, Cambridge 1939} (LFM; Diamond, 1975). Recent secondary literature includes e.g.\ M\"{u}hlh\"{o}lzer (2006, 2008, 2010, 2012), Schroeder (2020), Floyd (2021), Scheppers (2023),  and Bangu (2025). } which he extended to the philosophy of language.
In particular, \W's reflections on the analogy between mathematics and chess during his  middle period were pivotal in arriving at this late philosophy, as is especially clear from the following comment on Frege:\footnote{Secondary literature may be traced back from Lawrence (2023). Kienzler (1997) remains irreplaceable. }
  \begin{quote}
\begin{small}
Frege ridiculed the formalist conception of mathematics by saying that the formalists confused the unimportant thing, the sign, with the important, the meaning. Surely, one wishes to say, mathematics does not treat of dashes on a bit of paper. Frege's idea could be expressed thus: the propositions of mathematics, if they were just complexes of dashes, would be dead and utterly uninteresting, whereas they obviously have a kind of life. And the same, of course, could be said of any proposition: Without a sense, or without the thought, a proposition would be an utterly dead and trivial thing. And further it seems clear that no adding of inorganic signs can make the proposition live. And the conclusion which one draws from this is that what must be added to the dead signs in order to make a live proposition is something immaterial, with properties different from all mere signs.

But if we had to name anything which is the life of the sign, we should have to say that it was its \emph{use}.  (Wittgenstein, 1958, p.\ 4).
 \end{small}
 \end{quote}
 In other words, Frege (allegedly)  just saw two possibilities: either symbols refer to something in reality, in which case the game is meaningful (which, in his view,  doesn't apply to chess since it lacks ``thoughts'',  whose  absence supposedly blocks any analogy with mathematics), or they don't (which for Frege applies to chess but not to mathematics), in which case the game is meaningless. 
\W's point, then, is that Frege overlooked the possibility that even \emph{a priori} meaningless symbols might ``come alive'' by their use, as governed by the rules they are subject to.\footnote{As noted by Kienzler (1997), \S 4a, \W\ himself overlooks or ignores the fact that Frege (1903) \emph{does} note this `other possibility', for in footnote 1 on page 83 (which is part of \S71) he says: `Of course, there is also an opinion according to which numbers are neither symbols that mean something nor nonsensical meanings of such symbols, but rather figures that are handled according to certain rules, for example like chess pieces. According to this, the numbers are neither aids for research nor objects of observation, but rather objects of handling. This will have to be checked later.'  Indeed, in \S95 Frege (1903) 
complains that the rules of chess do not endow the chess pieces with any \emph{content} that would be the consequence of these rules, `like the name ``Sirius'' designates a certain fixed star.' This suggests a stubborn inability or refusal to see that the rules themselves comprise the meaning of chess (even though the pieces are meaningless); which was \W 's point. 
} 
 In sum:
\begin{itemize}
\item   For Frege, the use of symbols \emph{follows from their meaning}, given by their external referents; 
\item  For \W, the use of symbols, as determined by certain rules,  \emph{is their meaning}.
    \end{itemize}
\W\ also echoed our last quote from Hilbert to Frege, now in explicit reference to chess:
  \begin{quote}
\begin{small} 
It is, incidentally, very important that by merely looking at the little pieces of wood I cannot see whether they are pawns, bishops, castles, etc. I cannot say, `This is a pawn \emph{and} such-and-such rules hold
for this piece.' Rather, it is only the rules of the game that \emph{define} this
piece.
 A pawn \emph{is} the sum of the rules according to which it moves (a
square is a piece too), just as in language the rules of syntax define the
logical element of a word.
(Wittgenstein, 1967, p.\ 134).
 \end{small}
 \end{quote}

 During the 1930s \W\ moved from what has been called a  \emph{calculus conception} of mathematics, partly inspired by the analogy with chess, to a \emph{language-game conception}.\footnote{See Gerrard (1991), p.\ 127.}
This accompanied (and probably  even induced) a similar move in his philosophy of language,\footnote{ See Kienzler (1997) and Kuusela (2019) for his move in the philosophy of both mathematics and language.} which in fact is a more important source for us than his (unfinished) philosophy of mathematics.
Very briefly and with hindsight, in his \emph{Philosophical Investigations} he replaced \emph{referential} theories of meaning by what we reconstruct as \emph{inferential} ones, via the intermediate device of \emph{language games}. 

 Although \W\ himself refrains from giving a  definition of those--and would  surely regard any such definition as misguided, since games, languages, and language games are among his main examples of \emph{family resemblances}, which somehow defy definition--we  try it nonetheless:\footnote{\W\  introduced  language games as a  tool of his analysis of language in his middle period (notably in the \emph{Blue and Brown Books} from 1933--1935), generalizing this concept  in the \emph{Philosophische Untersuchungen}. The closest \W\ himself comes to at least a  characterization of language games is his list of examples in \S 23 of the PU. }
  \begin{enumerate}
\item A language game is a \emph{practice}  where certain words and symbols are \emph{used}.
\item The \emph{meaning} of (most) words and symbols is given by their \emph{use} within such a practice.
\item This \emph{use} is determined by specific \emph{rules} (forming the \emph{grammar} of the language game).\footnote{
 See especially \S\S 185--242 of the PI. For a brief discussion we  recommend  M\"{u}hlh\"{o}lzer, (2010), \S I.5.
Kuusela (2019), Chapter 6, holds that
language games \emph{need} not  be based on rules; but those in mathematics are.} 
\item These rules are \emph{inferential}: that is, the meaning of sentences lies in their inferential role.
\end{enumerate}
Though it seems compatible with (late) \W,  the last point was made more explicitly by
 Brandom (1995, 2001) concerning language and social practices as a whole.\footnote{This is the subject of an immense literature. Brandom (2007), Peregrin (2014) and Beran, Kolman, and Kore\v{n} (2018) are good points to start.  See also Haaparanta (2019) and Wischin (2019)  comparing  Brandom and \W.} 
Whatever the value of inferentialism in (natural) language, it does seem appropriate to logic (where the original inspiration came from),\footnote{The original sources of Brandom's inferentialism was Gentzen's proof system of \emph{Natural Deduction} in logic (see e.g.\  Von Plato, 2013),  where each logical symbol has an introduction rule and an elimination rule. These are seen as rules of inference for its use, from which its `usual' meaning is supposed to follow. This was in fact also suggested by \W, for example in his (1967), VII.30, as quoted below in our main text.
  Garson (2013),  Peregrin (2019), and Warren (2020) are inferentialist accounts of logic. }
and, though less easily implementable, to mathematics---though
surely, considerable work remains to be done in developing a full and satisfactory inferentialist account of pure mathematics.\footnote{Some steps were taken by Warren (2020) in support of his Conventionalism, but already our account of truth in \S\ref{truth}  differs from his, and also the role of (notably implicit) definitions as well as the origin of axioms, both \`{a} la Hilbert,  need to be clarified. 
  See our appendix below for further details and its relationship with Hilbert's formalism.  }  For the moment we take this possibility for granted and proceed. The question also arises \emph{which} language games in the above sense correspond to some form of mathematics; we let this question be answered by mathematical practice and ideas of family resemblance. The answer has to be fluid, if only because the rules of mathematics--and hence what came to be accepted as mathematics--have regularly  been subject to change; and despite  a century of apparent stability  they will surely change again (see also the end of \S\ref{truth} below). 
    
Furthermore, the key point found in both (late) \W\ and Brandom is that language--or mathematics--is primarily a (rule-governed) \emph{practice}; if one looks for a ``foundation'' of  language--or of mathematics--then this practice (rather than some logical formalism)
   \emph{is} its foundation.\footnote{For \W's philosophy of mathematics this view is developed in detail in  M\"{u}hlh\"{o}lzer (2010).}   
As already mentioned, such an  attitude towards mathematics seems to make it difficult to explain how mathematics can be such a powerful tool for describing the physical world.
Though we defer a detailed analysis to successor paper (Landsman, 2025ab), as a first step towards understanding  applications of at least basic arithmetic we again turn to \W. Among the coherent fragments of his philosophy of mathematics is the idea that
 mathematics does not provide \emph{representations} of ``reality'' (at least not primarily), but  \emph{yardsticks} to \emph{measure} 
  reality:
  \begin{quote}
\begin{small}What I want to say is: mathematics as such is always measure, not the thing measured.\footnote{'See also Kuusela (2019), \S4.4, for similar views on the philosophy of language during \W's middle period (notably in the \emph{Blue Book}).}
 (\W, 1969,  \S III.75h) 
   \end{small}
\end{quote}  
Conversely,  although mathematical results  typically \emph{originate} in experience, they are not empirical themselves, but should rather be seen as  `empirical propositions hardened into a rule':\footnote{See Schroeder (2020), Chapter 7, and Bangu (2025), Chapter 3, for detailed analysis of this passage.}
 \begin{quote}
\begin{small}
It is as if we had hardened the empirical proposition into a rule. And now we have, not an hypothesis that gets tested by experience, but a paradigm with which experience is compared and judged. And so a new kind of judgment.

For one judgment is: `He worked out $25 \times 25$, was attentive and conscientious in doing so and made it 615'; and another: `He worked out $25 \times 25$ and got 615 out instead of 625.'
But don't the two judgments come to the same thing in the end?

The arithmetical proposition is not the empirical proposition: `When I do \emph{this}, I get \emph{this}'--where the criterion for my doing \emph{this} is not supposed to be what results from it.
\\
(\W, 1969, \S VI.22bcd)
  \end{small}
\end{quote}  
Once the result has become a rule, it has become a piece of mathematics that as such is no longer subject to checks by experience: it is `put into the archives'.\footnote{See LFM,  Lecture IX, p.\ 107, further discussed in \S\ref{truth} below.}  This resolves the tension between mathematics as being timeless, non-spatial, acausal, etc.,  and yet applicable to our causal world in space and time.\footnote{The timeless of mathematics seems denied by Brouwer's  intuitionism, and perhaps some forms of constructivism.} This tension is due to a confusion between two different language games, namely pure mathematics, within which theorems are proved or looked up in some reliable book, and applied mathematics, in which theorems are compared with reality. As \W\ warned:
       \begin{quote}
\begin{small}
(\ldots) the sentence seems odd  only when one imagines it to belong to a different language-game from the one in which we actually use it. (\W, 2009, \S195). \end{small}
\end{quote}
 
 On the other hand, for all his earlier sympathy for the analogy between mathematics and chess and his emphasis on rules, later \W\ saw mathematics as applied \emph{by definition}:\footnote{On this topic see also Gerrard (1991),  M\"{u}hlh\"{o}lzer (2010), and Dawson (2014). } 
  \begin{quote}
\begin{small}
I want to say: it is essential to mathematics that its signs are also employed in \emph{mufti}.
It is the use outside mathematics, and so the meaning of the signs, that makes the sign-game into mathematics.
(\W, 1969, \S V.2)
 \end{small}
\end{quote}
This view is  problematic in modern mathematics, which is grounded in the autonomous development of mathematics from the 19th century onwards (see footnote \ref{Grayfn}). Likewise, his comments on the transition from the empirical origins of mathematics to a rule-governed activity reviewed above have an extremely limited scope in  modern mathematics and mathematical physics; and even within his elementary context \W\ had no acceptable theory of applied mathematics.

The  relationship between pure and applied mathematics, which is a central question for us, was  a major theme for Hilbert, who (unlike \W) was not just familiar with most or all of the mathematics and mathematical physics of his time; he had created or inspired much of it. 

Was Hilbert a ``formalist''?  He emphasized the importance of axiomatization throughout his career;  in the context of what is now called ``Hilbert's program'' this even came to include the rules of proof. As such, on a par with
Euclid he was a leading contributor  to the idea that mathematics is governed by (inferential) rules.\footnote{See  Mancosu, Zach, and Badesa (2009), and Ewald (2018).} 
But it  was one of his deepest conceptual insights that
\emph{both the rigour and the applicability of mathematics originate in axiomatization:}\footnote{Pasch, whose work foreshadowed Hilbert's in various ways, had  similar ideas (Pollard, 2010;  Schlimm, 2010).}
   \begin{quote}
\begin{small}
Mathematics has a two-fold task here: On the one hand, it is necessary to develop the systems of relations and examine their logical consequences, as happens in purely mathematical disciplines. This is the \emph{progressive task} of mathematics. On the other hand, it is important to give the theories formed on the basis of experience a firmer structure and a basis that is as simple as possible.

 For this it is necessary to clearly work out the prerequisites and to differentiate exactly what is an assumption and what is a logical conclusion. In this way, one gains clarity about all unconsciously made assumptions, and one recognizes the significance of the various assumptions, so that one can overlook what modifications will arise if one or the other of these assumptions has to be eliminated. This is the \emph{regressive task} of mathematics
 (Hilbert, 1919/1992, pp.\ 17--18). 
  \end{small}
\end{quote}
By axiomatization,  Hilbert meant the identification of certain sentences (becoming axioms) that form the foundation of a specific field in the sense that its theoretical structure (\emph{Fachwerk}) can be (re)constructed from the axioms via logical principles.  
The epistemological status of the axioms differs between various fields of mathematics. For example, Hilbert considered geometry initially a natural science that
emerged from the observation of nature (i.e.\ experience), which then turned into a mathematical science through axiomatization (Corry, 2004, p.\ 90).
This does not mean that he treated the axioms of geometry  as definitive, let alone as ``true" (as Euclid \emph{cum suis} had done).\footnote{See Pulte (2005) for a  history of the interpretation of mathematical axioms related to physics.}   Especially in physics Hilbert often stressed the tentative and malleable nature of axiom systems:  \begin{quote}
\begin{small}
As can be seen from what has been said so far, in physical theories the elimination of contradictions that arise will always have to be done by changing the choice of axioms and the difficulty lies in making the selection in such a way that all observed physical laws are logical consequences of the selected axioms.
  (Hilbert, 1918, p.\ 411)
 \end{small}
\end{quote}
For Hilbert, the axiomatization of physical theories is therefore never a static process: it moves on as physics itself moves on (Corry, 2004; Majer, 2014). Axiomatization may lead to the exposure of contradictions via a purely logical analysis, whose removal is then an important step forward. Indeed, as Majer powerfully summarized Hilbert's view on the axiomatization of physics:
 \begin{quote}
\begin{small}
physical theories live, as it were, on the border of inconsistency   (Majer, 2014, p.\ 72)
 \end{small}
\end{quote}
As Majer explains Hilbert, this is a consequence of an important difference between mathematics and physics in so far as axiomatization is concerned:
the former usually considers single disciplines in what he calls `maximal conceptual purity',  whereas
the latter often combines and intertwines a number of mathematical theories into a single highly complicated physical theory.

Axiomatization, then, contributes in two very different ways to the \emph{rigour} of mathematics:
\begin{enumerate}
\item  via syntactic proofs from the axioms (whose symbols remains uninterpreted);
\item via the axiomatization of sufficiently mature informal theories of mathematics.\footnote{Even in mathematics itself Hilbert acknowledged
the appearance of contradictions as a historical phenomenon. But unlike physics,  he apparently found contradictions unacceptable in mathematics, give his obsession of proving the consistency of classical mathematics. This marks a major difference with \W, whose cheerful acceptance of inconsistent theories, repeated comments on the indeterminateness of decimal expansions (e.g.\ of $\pi$),
 relaxed attitude towards the possibility of rejecting a correct proof, and whose insisting that proofs change or even define the nature of what was proved, sound  out of touch with modern mathematics and hence are inappropriate for our program.
}
\end{enumerate}
Similarly, axiomatization is also the key to the \emph{applicability} of mathematics, namely:
\begin{enumerate}[resume]
\item via the axiomatization of sufficiently mature theories of physics, space, quantity, etc.
\end{enumerate}
 In fact, it seems neither possible nor necessary to sharply distinguish between the second and third activities:
 for example, are Euclid's axioms (more precisely: his so-called postulates and common notions--whatever their clarity and worth from a modern point of view) attempts to axiomatize earlier informal geometry, or some physical theory of space?  Even the axiomatization of set theory in the early twentieth century brought rigour into both the  informal set theories of Riemann, Dedekind, and Cantor, and the genuine efforts by Frege, Russell,  and others to understand sets as ingredients of the physical universe or at least the human mind (Ferreir\'{o}s, 2008).
 
 \emph{Thus the key difference is between numbers 1 on the one hand and 2--3 on the other}:
the first is formal and focuses on proofs, whereas numbers 2 and 3  both take us outside (formal) mathematics.\footnote{ We follow Tait (1986) in seeing models in set theory as  internal to mathematics, and hence the distinction between syntax  and their interpretation in set theory  is  irrelevant for our theme. \label{tarski}} Thus it would reflect the spirit of Hilbert's `zweifache Aufgabe' (two-fold task) of mathematics to only list two sources of rigour in mathematics and the mathematical sciences:
\begin{enumerate}
\item[(i)] Defining mathematical theories by \emph{finding} appropriate axioms from heuristic considerations;
\item[(ii)] Proving theorems, \emph{given} these axioms (including deduction rules, themselves axiomatized). 
\end{enumerate}
\section{Mathematics as a rhododendron of language games}\label{MMLG} 
Despite his insights into the empirical sources of axioms and his impressive record in mathematical physics, even Hilbert hardly bridged the gap between his formalist views of pure mathematics (even if this was inspired by physics or other applications) and the mathematical description of natural phenomena: he relied on the  notion of  ``\emph{pre-established harmony}", a philosophical or even theological doctrine with roots in the monadology of Leibniz.\footnote{See Pyenson (1982),  Kragh (2015), and  Corry (2004).}   \W, on the other hand, left us both the impressive but far too narrow view that 
 mathematics consists of   `empirical propositions hardened into a rule' (see \S\ref{FHW}), and his famous description of what mathematics is:
  \begin{quote}
\begin{small}
Mathematics is a \textsc{motley} of techniques of proof. -- And upon this is based its manifold applicability and its importance.\footnote{Anscombe famously translated \W's `buntes Gemisch' as 
`motley', which is  the traditional costume of the court jester or fool.
 See M\"{u}hlh\"{o}lzer (2005), p.\ 66, footnote 15, for a critique of this translation.}
 (\W, 1969, \S III.46a)
 \end{small}
\end{quote}
 This juxtaposition of mathematical proof and applicability seems bizarre--unless one understands what \W\ means by `Beweistechniken' (proof techniques); the passage goes on as follows:
  \begin{quote}
\begin{small}
But that comes to the same thing as saying: if you had a system like that of Russell and produced systems like the differential calculus out of it by means of suitable definitions, you would be producing a new bit of mathematics.

Now surely one could simply say: if a man had invented calculating in the decimal system--that would have been a mathematical invention!--Even if he had already got Russell's Principia Mathematica.
 (\W, 1969, \S III.46bc)
 \end{small}
\end{quote}
In other words, by `Beweistechniken' \W\ means systems of definitions that, together with the logical deduction rules in \emph{Principia Mathematica} form the basis of new pieces of mathematics and hence (potentially) of new applications. If \W\  hadn't disliked set theory so much, the opening quote of this section could simply be that mathematics is a motley of its various branches, formalized within set theory. And this richness is indeed the key to its manifold applications and importance. Moreover, we may combine the first quote with a later one:
  \begin{quote}
\begin{small}
Logical inference is part of a language-game. (\ldots) We can conceive the rules of inference—I want to say—as giving the signs their meaning, because they are rules for the use of these signs
  (\W, 1969,  \S VII.30)
  \end{small}
\end{quote}
The suggestion then arises that mathematics is a motley of language games; although \W\ himself never seems to have claimed this in general,  the spirit of the idea is noticeable in both his \emph{Remarks on the Foundations of Mathematics} and the \emph{Philosophical Investigations}, and he does identify a few  special cases as such.
 But the setting in which he does so remains very limited. 

For a successful version  of this idea, one should incorporate all that mathematics involves: 
\begin{enumerate}
\item \emph{A long history:}  from numerical tables in Mesopotamia almost 4000 years ago to the rigorous concept of a function in the 19th and 20th centuries;
from quantitative methods of surveying to Riemannian geometry;  from counting to class field theory, \emph{et cetera}. It has thereby led to:
\item A number of different \emph{formal foundations of mathematics}, like ZF or ZFC or BNG set theory, intuitionistic set theory, $\lambda$-calculus, topos theory, homotopy type theory, \emph{et cetera}.
\item Within each of these foundational systems: a wide collection of mathematical theories (also called areas, branches, disciplines, or fields),\footnote{See the \emph{Mathematics Subject Classification} (MSC) at \url{https://mathscinet.ams.org/mathscinet/msc/msc2020.html}, or the `Branches of Mathematics' listed in
 Gowers (2008). 
} each with its own community, goals, and standards of proof. One may also think of Peano arithmetic or (Hilbert-style) Euclidean geometry.
These areas typically also overlap (e.g.\ Lie groups combine group theory and differential geometry; functional analysis combines linear algebra and topology, etc.). 
Following Hilbert (see \S\ref{FHW}) we find it hard to maintain the traditional distinction between  ``pure'' and ``applied'' mathematics (although many mathematics departments  do!).
\item Associated \emph{notions of proof} ranging from the informal reasoning of ancient Babylonian and Chinese mathematicians to the pseudo-axiomatic setting of Euclid (which lacked explicit rules of deduction) to the advanced logical apparatus of Frege, Russell, Hilbert, and G\"{o}del.
But even the logic differs not only between the formal foundational systems just mentioned (and others), but also includes considerable diversity in what is being tolerated within each of them, from informal rigour to  bending the rules.\footnote{ See M\"{u}hlh\"{o}lzer (2006) and Floyd (2023) for \W's notion of `surveyability' of a proof, which includes the criterion that its reproduction must be `an easy task'. This 
may be held against the  proof of  $1+1=2$  in \emph{Principia Mathematica} *54.43, which including all preparation takes hundreds of pages. Fully written Hilbert-style formal proofs of more complicated theorems would share this fate--though Dutilh Novaes (2011) identifies eight different ways in which formality and rules can be interpreted.
An informal proof is more likely to be understandable but is not rigorous, whereas a formal proof will hardly be understandable if only because of its length, which also  increases the probability of error (Avigad, 2021).
 In some cases  informal proofs are pointless, as in the four-colour theorem; we suggest that computer-assisted proofs and computer-verified proofs have their own language game (as defined below). 
}
See also \W's quote above.
\item The meta-theory of the axiomatized theories (i.e.\ Frege's ``theory of the game''), including \emph{both} formal aspects like proof theory \emph{and} informal aspects like ``the strategy of the game''.\footnote{Developing the formals aspect of this, i.e., metamathematics, was of course Hilbert's achievement.}
\item Applications of individual branches of mathematics to physics and other disciplines. 
\end{enumerate}
Expanding the notion of a game, we   answer our title question `\emph{Is mathematics like a game?}' by:
\begin{center}
\fbox{\emph{Mathematics is a rhododendron  of  language games of a very specific (formalized) kind.}}
\end{center}
Here our word `rhododendron' replaces \W's `motley' in order to emphasize that the structure of mathematics we propose has both multiple roots and various branches above these. The roots correspond to the various possible formal foundations, as just listed.\footnote{\W\ would already regard this starting point as misguided! Cf.\ M\"{u}hlh\"{o}lzer (2010) and Scheppers (2023).}  From each such root,  such as ZFC set theory, further pieces of formalized mathematics branch out, again as  listed.

Most of mathematics takes the form of three  language games played on \emph{the same}  theories:
\begin{enumerate}
\item The cleanest mathematical language game is the development of theorems and proofs---where different proof systems may be used, defining different language games (cf.\ \S\ref{truth}). 
\item The second one is Hilbert-style metamathematics, yet with one version for each root.\footnote{If (against our advice) mathematics is seen as a game, then metamathematics is  Frege's ``theory of the game''. Whatever one's philosophy of mathematics, Frege was right that the ``theory of the game'' cannot be about meaningless symbols since it is about mathematical proofs;  and this interprets the symbols. Nonetheless, metamathematics is determined by inferential rules (like mathematics itself) and the corresponding ``meaning is use'' semantics returns this interpretation. Hence the  metamathematical language games share this semantics with all the other mathematical language games and in that sense they are all on the same par. A similar comment applies to the theory of any game with a mathematical structure (which includes
 the metamathematics of some given root of the rhododendron).}
\item The third is applied mathematics in the specific form summarized in the Introduction.
\end{enumerate}
The first is ``mathematicians's mathematics''. The second is  played by logicians and philosophers. The third, which practiced  by applied mathematicians and mathematical physicists,  is obviously restricted to `applicable'' branches of mathematics (which of course expand in the course of time).

Hilbert's formalist emphasis on the meaninglessness of mathematical symbols applies to the first, but only \emph{in so far as proofs and other formal aspects of axiom systems are concerned} (such as consistency and completeness). 
The analogy between mathematics and chess applies here. This analogy is not uniquely definable, but an attractive version is the one  proposed by Weyl (1926):\footnote{ The last point, which is the idea of implicit definition,  was not  mentioned by Weyl (1926) and should be attributed to Hilbert and others; see footnote \ref{IDfn}.  What is admittedly missing in the analogy between mathematics and chess is a translation of the goal of \emph{winning} in chess: there seems to be no analogue of checkmate in mathematics (although there is an emotional analogue of resigning, i.e., ``giving up'', after repeated failure to prove some theorem). Indeed, in the latter the goal is to establish the counterpart not of a winning position but of an arbitrary legal position (i.e.\ a theorem). Perhaps the shared aspect of beauty in both games of chess and proofs  somewhat compensates for this discrepancy. }
\begin{itemize}
\item  The \emph{axioms} of some theory are analogous to the starting position of a game of chess;
\item The  \emph{deduction rules} (\`{a} la Natural Deduction) are analogous to the possible moves;\footnote{
What we have in mind here is that  the axioms are supposed to describe some specific mathematical theory (such as set theory, or arithmetic, or Euclidean geometry) whereas all deduction rules are logical in character and, perhaps with a few exceptions, are universal for all fields of mathematics (like Euclid's common notions). See e.g.\ von Plato (2017). In contrast,
in a Hilbert-style calculus  (Hilbert \& Ackermann, 1928)  \emph{modus ponens} is the only deduction rule
whilst the other deduction rules in Natural Deduction are seen as axioms. This calculus does not fit our metaphor. }
\item A \emph{sentence} (as defined in logic)  is analogous to \emph{some} position on a chess board;
\item  A \emph{theorem} is like a \emph{legal} position in a correctly played chess game;
\item A \emph{proof} is like a game leading to that position, played according to the rules;
\item A \emph{definition} resembles the idea that chess pieces are defined by the rules of chess.
\end{itemize}
 Given the formal notion of proof developed by by Frege, Russell, and Hilbert, in which (unlike in Euclid--let alone 17th and 18th century mathematics) not only the axioms but also the rules of deduction are formalized and  stated, the rules of this language  game are clearly inferential.

 The second language game  is one level above the previous one but it also squarely lies on the formal and inferential side: although  the object of investigation is a mathematical theory, its symbols remain uninterpreted. It is played on the same theories as the previous one (for example, arithmetic, as in ``Hilbert's program''), but it need not follow the same logic as the original game.\footnote{Here one need not think of the full scope of Hilbert's program; G\"{o}del's completeness theorem is already metamathematical, as is its special case for propositional logic (Zach, 1999). }

The third one, which is crucial in understanding  the relationship between  mathematics and the physical world, is based on rules that are inferential in a much less obvious way, for which we again refer to the Introduction for a summary and to Landsman (2025ab) for a full development. 

And of course there are numerous other language games  that are especially relevant to mathematical practice; for example, in either trying to find proofs or in creating new mathematics even the most stubborn formalist will look for interpretations and perhaps visualizations in both applied and non-applied (or not-yet-applied) fields, for links with different fields of mathematics, etc. Similarly, learning mathematics is a language game by itself (much as learning a language is, as analyzed in the early parts of the \emph{Philosophical Investigations}). And so on and so forth.
 \section{Truth}\label{truth}
  \begin{center} 
\begin{small}
 The truth predicate then preserves his contact with the world, where his heart is. (Quine, 1986, p.\ 35)
 \end{small}
\end{center}
What does the ``pure'' or ``formalist'' language game imply for the concept of \emph{truth} in mathematics?
  On the one hand this is a difficult question, since  we see mathematics as a human practice, for whose rules  there were always many different possibilities and choices, even when these rules were inspired by empirical phenomena. On the other hand, mathematical practice suggests that ``truth'' is merely a \emph{fa\c{c}on de parler}, in that for most ``working mathematicians'' ``$p$ is true'' simply means that $p$ is a theorem.
Anything beyond this meets the scathing comment of Bourbaki:
\begin{quote}
\begin{small}
Mathematicians have always been sure that they prove ``truths'' or ``true propositions''; such a conviction can obviously only be sentimental or metaphysical. (Bourbaki, 1994, p.\ 11)
\end{small}
\end{quote}
We agree, but some place for ``truth'' remains in mathematics; it just needs to be relocated. 

Let us return to chess for inspiration. It seems meaningless to say that a position $p$ in chess is ``true''. But it does make sense to claim that $p$ is \emph{legal}, in that it arose from a  game played according to the rules $R$. This claim, call it $R\vdash p$, rather than $p$ itself,  could be said to be true or false, and this can be established by a proof in the form of an actual (legal) chess game leading to $p$, cf.\ \S\ref{MMLG}. 
 In normal games $p$ even arises in this way; in so-called retrogade chess problems one has to reconstruct $p$.
Similarly,  in our non-referential ideology mathematical theorems cannot be true either, since there is no objective state of affairs they could describe correctly.\footnote{We repeat the point already made in footnote \ref{tarski}: Tarski's concept of truth as defined in model theory is internal to pure mathematics and has little or nothing to do with the notion of truth sought by the Platonists or naturalists. It will play a minor role in the discussion of G\"{o}del's theorems below.}
Moreover, the kind of truth conventionalists aspire to is, in our view, covered much better by our proposal below than by declaring theorems of propositions themselves to be true.
Like in chess, truth in mathematics cannot lie in sentences $\phv$ (such as closed formulae in first-order logic), but only in
 claims $T\vdash\phv$ stating \emph{that $\phv$  is a theorem within an ambient theory $T$} (which is supposed to include rules of inference). 
  And this is the case (by definition) iff there exists a proof of $\phv$ according to the rules of $T$.
 Thus the only thing we can say about mathematical truth in our framework is this:
    \begin{center}
\fbox{\emph{Mathematical truth resides not \emph{in} theorems but in claims \emph{that} some sentence is a theorem. }}
\end{center}   
 This makes a proof of $\phv$ in $T$ the truth-maker of the truth-bearer $T\vdash\phv$. Our only compromise towards realism  is our belief that such truth (or falsehood) is a matter of \emph{fact}, whether or not it is \emph{known}.\footnote{This might be varied by defining  $T\vdash\phv$ to be true if a proof of $\phv$ \emph{is} known, as in intuitionistic mathematics.}  But this is not the truth of platonism (or of  naturalistic views of mathematics),
  which concerns $\phv$ rather than $T\vdash\phv$, backed by 
 a correspondence theory of truth. We also reject other approaches which argue 
  that a sentence $\phv$ \emph{itself} (as opposed to $T\vdash\phv$) is true iff $\phv$ has a proof.\footnote{See e.g.\ Dieudonn\'{e} (1971) and Tait (1986). See also  Appendix \ref{formalism} for truth in formalism and deductivism.
  The so-called BHK (Brouwer--Heyting--Kolmogorov) interpretation (or semantics) of intuitionistic logic is also often taken to mean that a sentence is true iff it has a proof (Artemov \& Fitting, 2021). We reject this, too, but even so one may still support the more modest BHK interpretation of the logical connectives in terms of proofs (van Atten, 2023). }  

 First, unless one believes that there is a single ``true'' foundational system for mathematics (such as ZFC set theory with additional cardinality axioms, as proposed by G\"{o}del),  such proposals endorse a \emph{coherence theory of truth} (Young, 2018), in which each such system would come with its own set of truths. As explained in \S\ref{MMLG}, we reject this (cf.\ Landsman, 2025ab). On our proposal, although people may differ about the virtues of different
foundational systems, given unambiguous concepts of inference and proof  they cannot rationally differ about the theorems in each of these.

Second,\footnote{See also Paseau \&  Pregel (2023), \S 9, and references therein. }
  according to  G\"{o}del's first incompleteness theorem,  for any $T$ (satisfying the usual assumptions) there are sentences $\phv$ such that neither $\phv$ nor $\neg\phv$  is provable in $T$. Yet in classical logic $\phv\vee\neg\phv$ is \emph{provable} for every  $\phv$. If this implies that $\phv\vee\neg\phv$ is \emph{true}, then for undecidable $\phv$ this would be the case without either $\phv$ or $\neg\phv$ being true, which is awkward.\footnote{This problem obviously  does not arise in intuitionistic logic, which G\"{o}del (1931) actually incorporated.} But since the claim `$T\vdash(\phv\vee \neg\phv)$' is clearly different from 
`$(T\vdash\phv) \mbox{ or } T\vdash(\neg\phv)$',  there is no argument to conclude from the truth of the former that the latter is true, and so even on  TND and the everyday understanding of `or' we are not forced to (wrongly) conclude that  either $T\vdash \phv$ or $T\vdash\neg \phv$ is true, 

Thus the possibility of assigning truth  to $\phv$  is challenged by G\"{o}del's incomplete\-ness theorems, whereas $T\vdash \phv$ faces no such problems.\footnote{As brought to our attention by Barteld Kooi, 
incompleteness does lead us to an asymmetry between truth and falsehood the naive approach (in which $\phv$ itself is true or false) does not have. Namely, if  $T\vdash\phv$ is true then it has
a \emph{truth}-maker in the form of a proof of $\phv$ from $T$; but if $T\vdash\phv$ is false, i.e., not true,  then there is
a \emph{false}-maker for it just in case $\phv$ is decidable in $T$ (in which case $T\vdash\neg\phv$ and hence 
 the false-maker is a proof of $\neg\phv$ from $T$). 
But this seems a lesser evil (for us) than some assumption of mathematical realism. } If Dummett's (1963/1978) famous (but controversial) argument for the `vagueness' of the concept natural number is correct,\footnote{See e.g.\ Engler (2025) and references therein.
See also Parsons (1990) for arguments similar to Dummett's. }  then the stance of arithmetical determinacy mentioned in the Introduction (and in its wake the more general `fundamental commitment of mathematics') is also \emph{weakened} by  G\"{o}del's incomplete\-ness theorems, although these are usually seen as a \emph{threat} to philosophies  like ours.\footnote{ See Weir (2010), \S4.III,
for formalism;  Warren (2020), \S11.VII, for conventionalism; and  Paseau and Pregel (2023), \S9, for
deductivism. See also \S\ref{formalism} below for our relationship to these philosophies or attitudes.} The following preamble is uncontroversial.
 Take $T = PA$ (i.e., Peano Arithmetic with first-order logic) to be specific, and let  $G_T$ be its G\"{o}del sentence, which  expresses its own unprovability in $T$.\footnote{In very naive discussions this sentence is  claimed to be true full stop, and since it cannot be proved (in $T$) this is supposedly an argument against formalism and for platonims and/or the superiority of the brain over any formal system. See Franz\'{e}n  (2005) for a critical discussion of this and  many other misunderstandings of  G\"{o}del's theorems.} Then it is well known that $G_T$ is undecidable (i.e., neither $T\vdash G_T$ nor $T\vdash\neg G_T$) \emph{if and only if $T$ is consistent}.\footnote{See e.g.\ Franz\'{e}n (2005) or Raatikainen (2022).}  Then:
 \begin{itemize}
\item  Arithmetical determinists take  the ``existence'' of the natural numbers $\mathbb{N}$ and their satisfaction of the PA axioms as given, whence PA is consistent. From this,  they are entitled to conclude that  the interpretation $[[G_T]]_{\N}$ of $G_T$ in $\N$ is true in their absolute sense (cf.\ \S\ref{Intro}). 
\item Without the commitment to arithmetical determinacy, 
$[[G_T]]_{\N}$ is merely true  in the formal sense of Tarski for model theory, where
 $\N$ is  a  construction within ZFC set theory.\footnote{Or in some fragment thereof in which the construction of $\N$ can be carried out and in  which the consistency of $PA$ can be proved. This only makes sense if 
ZFC (or the  fragment just alluded to) is consistent, which of course is a big if.} 
\end{itemize}
Thus the claim that $[[G_T]]_{\N}$ or $G_T$ is true in an ``absolute'' sense, i.e.,   beyond derivability
in some axiomatized theory (such as ZFC in the second case), must already \emph{assume} arithmetical determinacy:  the incompleteness theorems cannot be used to \emph{derive} it or even argue for it.\footnote{As does for example Connes in  support of his Platonism (Connes, Lichnerowicz, and  Sch\"{u}tzenberger, 2001). }Indeed: \begin{quote}
\begin{small}
The [arithmetical determinist] however,
operates with the notion of a model as if it were something that could be
given to us independently of any description: as a kind of intuitive conception which we can survey in its entirety in our mind's eye, even though we
can find no description which determines it uniquely. This has nothing to
do with the concept of a model as that concept is legitimately used in
mathematics. There is no way in which we can be `given' a model save by
being given a description of that model. If we cannot be given a complete
characterisation of a model for number theory, then there is not any other
way in which, in the absence of such a complete description, we could
nevertheless somehow gain a complete conception of its structure. (Dummett, 1978, p.\ 191)
\end{small}
\end{quote}
 This confirms the circularity in arguing that  G\"{o}del's theorems  enforce  arithmetical determinacy.
 Our account of this situation is as follows, \emph{assuming} PA is consistent. No sentence in PA = T, including $G_T$, has anything like a truth value. The claims $PA\vdash G_T$ and $PA\vdash \neg G_T$ are both false in our sense (i.e., not true). The ``truth'' of $[[G_T]]_{\N}$ is just the truth of $\N\vDash G_T$, seen as a theorem in ZFC or in some  weaker system, such as the proof system one obtains by adding the so-called $\omega$-rule to PA.\footnote{This rule, which goes back to Hilbert, states that $\phv(n)$ for all $n$ defined in PA as $0=0$, $1=S(0)$, $2=S(S(0))$, etc., 
 implies $\forall_x\phv(x)$. This rule  allows one to prove all true statements $\forall_x \phv(x)$ provided $PA\vdash \phv(n)$ for each $n=0, 1, 2 , \ldots$,
 at the expense of using an infinite number of assumptions. As such it distinguishes the standard model $\N$ of PA from all other (i.e.\ non-standard) models, in that 
 the $\omega$-rule holds for all $\phv(x)$ precisely in $\N$. The G\"{o}del sentence $G_T$ is of the form $G_T=\forall_x \phv(x)$, where each $\phv(n)$ for $n=0, 1, 2 , \ldots$  is a theorem of PA. The gap between
 the inability to prove $G_T$ in PA and the ability to prove it in $\N$ is therefore precisely bridged by the  $\omega$-rule.
 See e.g.\ Warren (2020), \S10.VII. \label{omegafn}} In view of our definition of truth we here effectively replace ``truth-talk'' (in the usual sense) by ``proof-talk'', with the crucial feature that by changing the proof system in passing from $PA\vdash G_T$ (which is false) to $\N\vDash G_T$ (which is true) \emph{we switched to a different language game}.\footnote{See also Kolman (2014, 2016), which makes a similar point in a different context. }

Still assuming consistency of PA, there is a non-standard model $\mathbb{N}'$ of PA in which  the interpretation $[[G_T]]_{\N'}$ of $G_T$ is false.\footnote{Continuing footnote \ref{omegafn}: each sentence $\phv(n)$ appearing in $G_T=\forall_x \phv(x)$ is true on $\N'$, yet $G_T$ is false in $\N'$.} This falsehood indicates that the concept of a natural number  is not sufficiently captured by PA,  and since this argument is independent of the choice of PA  arithmetical determinists  must agree that `we have a certain, quite definite, concept, which cannot be fully characterised just by the fact that we make certain assertions about it' 
and that  we `cannot characterise completely the meaning of
``natural number'' by specifying which arithmetical statements we are prepared to assert and which forms of inference within arithmetic we are prepared to accept.'\footnote{Dummett (1978), p.\ 186, 187.}
\emph{We} conclude that since PA--or some similar system, facing similar problems---is the only intuition we have about natural numbers, one cannot possibly claim that the natural numbers are ``defined'' or ``exist'' in some absolute sense.\footnote{We  recall the sad fact that Frege's life work of defining the natural numbers  failed even for the number one.} But \emph{arithmetical determinists}, instead of giving up their position,  conclude from this that at least in this case  \emph{meaning} (namely the absolute concept of natural numbers they have in mind) cannot be given by \emph{use} (according to the axioms and rules of inference). It is this way out that the controversial remainder of Dummett's argument tries to undermine,  so that Dummett's higher goal lies in defending a ``meaning = use'' semantics. 
 
 In any case, one surely needs to get used to the idea that say $7+5=12$ is neither true nor false (it is just not the kind of mathematical statement that has a truth value),\footnote{An anonymous referee highlighted the radical nature of our concept of truth by mentioning 
 the Sylow theorems for finite groups. Here the temptation to relate a purely mathematical claim to real things like apples seems absent, but this \emph{weakens} the pull to attach any truth label to such theorems. Like all mathematical objects, finite groups and their properties have a `permanence of identity through space and history' (see Rota quoted in the Introduction), but like numbers this identity is  given by specific definitions and other rules, as opposed to things that actually have the said properties. This \emph{increases}  the pull in the opposite direction of believing that theorems are statements about rules.
 }
  whereas the superficially similar but technically and conceptually very different claim $\mathrm{PA} \vdash (7+5=12)$,
  stating that $7+5=12$ is a theorem of Peano arithmetic, or equivalently that  $7+5=12$ is true \emph{in PA}, 
 \emph{is} true.\footnote{In some crazy theory $T$ where
$T\vdash (7+5=10)$, this would still be true on our criterion (as long as the proof in $T$ is correct!). 
This theory might be an interesting game but it would be a useless yardstick in applied mathematics. }

Don't seven apples add up with five apples to yield twelve apples? They do. But this expresses neither
$7+5=12$ nor $\mathrm{PA} \vdash (7+5=12)$: the former  is \emph{held against the apples  as a yardstick}, justified by the latter. This yardstick even seems to yields perfect results (which is rare and may be restricted to counting and elementary arithmetic), but already Aristotle  realized how much is involved in this: one must regard each apple as a unit, which deliberately  overlooks  that firstly each apple is divisible, and secondly that all apples are different. Even granting these idealizations, what we have is a match between empirical data and some mathematical theorem, viz.\ $7+5=12$.
Following \W, the latter is an `empirical proposition hardened into a rule'. But this rule  cannot inherit any kind of truth from the empirical propositions that originally inspired it (such as the counting of objects like apples), since that would confuse the physical world with the role of mathematics as a set of yardsticks invented by humans to understand it:
 \begin{quote}
\begin{small}
I am trying to show in a very general way how the misunderstanding of supposing a mathematical proposition to be like an experiential proposition leads to the misunderstanding of supposing that a mathematical proposition is about scratches on the blackboard.

Take ``$20+15=35$''. We say this is about numbers. Now is it about the symbols, the scratches? That is absurd. It couldn't be called a statement or proposition about them; if we have to say that it is a so-and-so about them, we could say that it is a rule or convention about them.--One might say, ``Could it not be a statement about how people use symbols?'' I should reply that that is not in fact how it is used--any more than as a declaration of love.
(Wittgenstein, LFM, Lecture XII, p.\ 112)
\end{small}
\end{quote}

The unity of mathematics emphasized by Hilbert provides an additional argument for the lack of truth of $7+5=12$. If this theorem were true, then every theorem in mathematics should be  true. This leads to a problem discussed earlier:  since different foundational systems may yield contradictory results,   just one of these  systems could be ``true''. The history of mathematics suggests this is dubious. Even if only one of them ultimately comes out be correct (e.g.\ since the others unexpectedly are  inconsistent), putting esoteric result in ZFC set theory about inaccessible cardinals on a par with  $7+5=12$ as both being ``true'' sounds equally wrong. 
The only way to get around these problems seems to be to treat all theorems from all foundational systems on a par;
but instead of declaring them all true, the ensuing notion of truth is expressed much better by saying that the claim $T\vdash \phv$ that $\phv$ can be deduced from $T$ is true, rather than $\phv$ itself. 

We see that even the simplest  theorems in arithmetic, involving very small integers, are no threat to
our notion of truth.  Our case against realism is even stronger if   large integers are involved; and still stronger for Euclidean geometry; and  stronger again if we use the advanced theories of mathematics  physics, about which only very naive physicist would say  their theorems are ``true''; see Landsman (2025ab) for further analysis. But our argument is uniform for all these cases.

 \W's (disquotational) concept of truth differs  from ours:  in \S 136 of the \emph{Philosophical Investigations} he identifies
 \emph{`$p$' is true}  with $p$ itself (and \emph{`$p$' is false} with not-$p$).
 But in a marked difference with the `fundamental commitment of mathematics'  he then adds that \emph{this concept of truth belongs to the rules of the language game in question}.\footnote{Similarly in \S6 of Appendix II of \W\ (1969). } 
Similarly, \W\ wrote:
 \begin{quote}
\begin{small} 
Mathematical truth isn't established by their all agreeing that it's true--as if they were witnesses to it. \emph{Because} they all agree in what they do, we lay it down as a rule, and put it in the archives.  (\W\ in Diamond, 1975,  Lecture IX, p.\ 107) 

Disputes do not break out (among mathematicians, say) over the question whether a rule has been obeyed or not. People don't come to blows over it, for example. (\W, 2009,  \S 240)
\end{small}
\end{quote}
\section{Conclusion}\label{Conc} 
 \begin{quote}
\begin{small} 
What is the conclusion to which we come? Modern mathematics and physics may seem to move in thin air. But they rest on a quite manifest and familiar foundation, namely the concrete existence of man in his world. (Weyl, 2009, p.\ 188) 
\end{small}
\end{quote}
Our  aim was to re-examine the question to what extent mathematics may be compared with a game. Trying to answer this question led us essentially to Weyl's conclusion just quoted.  To get there,  we combined certain  insights by Hilbert and \W\ (as ``inferentialised''  by Brandom). From Hilbert, we  took the idea that axiomatization enables both pure and applied mathematics (including mathematical physics). From \W, our main lessons is that   pure and applied mathematics correspond to different language games, both of which are non-referential. 
 The ``applied'' one also relies on his remarkable idea of using mathematical theorems or theories as yardsticks. 
 
 Mathematics also incorporates various other language games, together forming a structure we like to call a ``rhododendron'' (rather than \W's  ``motley'', which seems too flat). Various mainstream philosophies of mathematics find a peaceful place within this structure, except Platonism: Platonists and other realists will like little of our proposal, since they already reject our (and \W's)  starting point of ultimately grounding mathematics in human practice. 
  Formalists may have more sympathy for our view, since we not only defend what should be the deductivist concept of truth, but also moved the problem of understanding applied mathematics within a formalist framework a step forward.  Logicism, formalism, and intuitionism may no longer exist in their original form, but what is left of them is also welcomed within the rhododendron. 
We take this peaceful coexistence to be a major advantage of our proposal.
We propose this rhododendron of language games as an object of comparison itself: we invite readers to compare the mathematics they have in mind with this picture, to see where it agrees and where it deviates. 

Finally, although our paper is by no means meant as an analysis of \W's philosophy of mathematics (nor of Hilbert's), we hope to have implicitly answered certain objections to it  by firstly doing some cherry-picking (i.e.\ tacitly removing some of his more extreme and outdated views) and secondly integrating these marbles with modern ideas appropriate to contemporary mathematics (especially ideas of Hilbert's, hopefully without falling into \emph{his} traps either).
\appendix
\section{Formalism, deductivism, conventionalism, and pluralism}
In tis appendix we briefly compare our approach with certain interpretations (or even philosophies) of mathematics that are related to ours and  of which we have tried to incorporate some parts. 
 \subsection{Formalism and deductivism}\label{formalism}
Though often associated with Hilbert,\footnote{Indeed, the exposition of formalism by von Neumann (1931) is indeed entirely devoted to Hilbert's program.} there isn't a canonical notion of ``formalism''. Here are two caricatures from leading textbooks in the  philosophy of mathematics:
\begin{quote}
\begin{small}
The various philosophies that go by the name of `formalism' pursue a claim that the \emph{essence} of mathematics is the manipulation of characters. A list of the characters and allowed rules all but exhausts what there is to say about a given branch of mathematics. (Shapiro, 2000, p.\ 140). 

\emph{Formalism} is the view that mathematics has no need for
semantic notions, or at least none that cannot be reduced to
syntactic ones. \hfill (Linneb\o, 2017, p.\ 39)
\end{small}
\end{quote}
This is the version attacked by Frege; cf.\ his correspondence with Hilbert quoted in \S\ref{FHW}. But Hilbert  is a straw man;  for him it is only \emph{in the context of proofs and in the analysis of consistency of axiom systems etc.} that mathematics is  a deductive enterprise in which symbols have no meaning (outside the rules they are subject to). A  broader view of formalism is described by Detlefsen (2005).
Our own approach incorporates formalism on the side of  pure mathematics, but tries to balance it via Hilbert's emphasis on the informal meaning of symbols inherited from  the heuristic theories of either physics or mathematics inspiring most axiomatizations.  His implicit definitions, or \W's corresponding ideas about the meaning of symboles given by their use (`the life of the sign')  also gives meaning (or even `life') to formalism in a way that many philosophical discussions seem to overlook.\footnote{In addition, Hilbert's program of giving a finitist consistency proof of classical mathematics re-introduced content even in a purely formal setting, cf.\ Hilbert (1918). This is taken into account by the second language game in \S\ref{MMLG}. }
 We are not ``formalists'' if this is seen, as it historically has, in opposition to logicism or intuitionism; appropriately formalized (!) at least the latter is a valid mathematical language game on a par with classical mathematics as used by most formalists. 

Our concept of truth within the formalist language game expounded in  \S\ref{truth} is compatible with a philosophy of mathematics called \emph{deductivism}, famously summarized by Russell (1903) as: 
\begin{quote}
\begin{small}
Pure mathematics is the class of all propositions of the form ``$p$ implies $q$''.
\end{small}
\end{quote}
See  Paseau and Pregel (2023).
 Different ``deductivists'' had different concepts of truth. For example, Russell  saw the truth of logic and mathematics in the universe. Around 1900, Hilbert defined mathematical truth as flowing from axioms, which by themselves were deemed ``true'' if they were consistent (Paseau and Pregel, 2023, \S3.1).  Curry (1951)'s  Chapter II is called `The problem of mathematical truth' and starts with the statement that `The central problem in the philosophy of mathematics is the definition of mathematical truth'. It may then come as a slight disappointment that he simply
 puts truth in the theorems themselves, justified by a verification procedure that comes down to checking their proofs.  Lolli (1998) states that only \emph{logical} truth  matters, defined  as `truth under any interpretation whatsoever' which in turn he takes to mean `true under any notion of truth' (p.\ 118), adding that `the great success of mathematical logic is to have shown that all of logic is independent of a definition of truth. (p.\ 119). 
For  Weir (2010),  truth comes from the correctness of utterances, which in mathematics is guaranteed by provability. But unlike in our analysis it is still the utterance, i.e., the theorem, which is deemed true. \emph{Et cetera}.
\subsection{Pluralism and Conventionalism}
Our proposal is meant to give room to various ``philosophies'' or ``foundations'' of mathematics that are often seen to be mutually incompatible, such as classical mathematics (based on set theory) and the  formalism often associated with it, intuitionism, constructivism, etc. We see these as different language games. As such, it is clearly ``pluralist'' in character. But what does this mean?

Friend (2014) ends her monograph on mathematical pluralism with the following manifesto:
\begin{quote}\begin{small}
 One can be pluralist in different respects, at different levels and one's pluralism can be governed by different logical inclinations or hypotheses. In general, the pluralist aspires to the following virtues: unprejudiced observation of mathematical practice and a desire to encompass and accommodate as wide a variety of practices as is coherently possible. The inverse of these virtues are manifested when we insist on unique, simple, teleologically satisfying answers, beyond what the evidence will support. (\ldots) The pluralist position is meant to give a philosophical theory to support what is already happening in the philosophy of mathematics. Thus the position is `new' in the sense of not having yet been expressed this way in print, but it is `old' in the sense of being already implemented and understood, at some level, by some mathematicians and philosophers of mathematics. (Friend, 2014, pp.\ 241--242) 
\end{small} \end{quote}

In a more recent text,  Priest (2024) argues for the coexistence of various foundations (or `arenas') of mathematics (including even substructural logics), draws the analogy with games,  and even ends with the same quote of \W\ that appears in the opening of our \S\ref{MMLG}.\footnote{The first version  of our article, still available as
\url{https://arxiv.org/abs/2311.12478v1},  already contained this quote  before the publication of Priest (2024), which we came to know only in 2025.} 
Finally, a recent paper by Zalta (2024) opens with:
\begin{quote}\begin{small}
Mathematical pluralism can take one of three forms:
(1) every consistent mathematical theory consists of
truths about its own domain of individuals and relations;
(2) every mathematical theory, consistent or inconsistent, consists of truths about its own (possibly
unintersting) domain of individuals and relations; and (3) the
principal philosophies of mathematics are each based
upon an insight or truth about the nature of mathematics that can be validated. 
\end{small} \end{quote}

Our paper is  written in the same spirit, contributing to mathematical pluralism. Priest (2024) is closest to us, but none of these works  present a similar the analysis and  architecture of mathematics (based on Hilbert and \W) or has a comparable approach to applied mathematics. 

\emph{Coventionalism} is a philosophical position that is closely related to mathematical pluralism, as well as to our work. It has its roots in Poincar\'{e}, early Wittgenstein, and logical positivism (Ben-Menahem, 2006).  Given our analysis in \S\ref{FHW}, it should be no surprise that the species of conventionalism of interest to us is \empty{inferentialism}, in which the relevant conventions are (syntactic) inference rules. This is developed for natural language by Brandom and followers (see \S\ref{FHW} and references therein), whereas for logic and mathematics Warren (2020) is a major source; see also Garson (2013) for logic.\footnote{ Despite these books (and references therein), much work remains to be done in order to derive an inferentialist account of pure and applied mathematics; this will be taken up in the future.}
 All of this is  also indebted to late \W, as well as to Hilbert, and we largely align with it.  As far is truth is concerned, there seem to be two main directions, of which we follow Ben-Menahem (2006) in thinking of conventions as hypothetical conditions rather than as freely postulated truths, as in Warren (2020).
 See Landsman (2025ab) for details. 
\addcontentsline{toc}{section}{References}
\begin{small}

\end{small}
\end{document}